\documentclass[12pt]{amsart}
\usepackage{amssymb,amsmath,latexsym, enumerate}

\usepackage{color}

\numberwithin{equation}{section}
\theoremstyle{plain}
\newtheorem{thm}{Theorem}[section]

\newtheorem{ex}[thm]{Example}
\newtheorem{cor}[thm]{Corollary}
\newtheorem{prop}[thm]{Proposition}

\newtheorem{df}[thm]{Definition}

\begin{document}
\newcommand{\hgt}{\operatorname{ht}}
\newcommand{\Hht}{\operatorname{Hht}}
\newcommand{\Ext}{\operatorname{Ext}}
\newcommand{\Hom}{\operatorname{Hom}}
\newcommand{\Supp}{\operatorname{Supp}_R}
\newcommand{\Ass}{\operatorname{Ass}_R}
\newcommand{\Att}{\operatorname{Att}_R}
\newcommand{\Min}{\operatorname{Min}_R}
\newcommand{\depth}{\operatorname{depth}}
\newcommand{\pd}{\operatorname{pd}}
\newcommand{\ann}{\operatorname{ann}}
\newcommand{\grade}{\operatorname{grade}}
\newcommand{\Grade}{\operatorname{Grade}}
\newcommand{\Spec}{\operatorname{Spec}}
\newcommand{\lc}{\operatorname{lc}}

\title[Finite Minimal Components]{On rings for which finitely generated ideals have only finitely many minimal components}
\author{Thomas Marley}
\address{Department of Mathematics\\
University of Nebraska-Lincoln\\
Lincoln, NE 68588-0130} \email{tmarley1@math.unl.edu}
\urladdr{http://www.math.unl.edu/\textasciitilde tmarley1}

\date{December 20, 2005}

\begin{abstract}  For a commutative ring $R$ we investigate the property
that the sets of minimal primes of finitely generated ideals of
$R$ are always finite.  We prove this property passes to
polynomial ring extensions (in an arbitrary number of variables)
over $R$ as well as to $R$-algebras which are finitely presented
as $R$-modules.
\end{abstract}

\maketitle

\section{Introduction}

In \cite{OP} Ohm and Pendleton examine several topological
properties which may be possessed by the prime spectrum of a
commutative ring $R$. One of these properties, denoted in \cite{OP}
by {\it FC} for `finite components', is that every closed subset of
$\Spec R$ has a finite number of irreducible components, or
equivalently, that every quotient $R/I$ has a finite number of
minimal primes.  In this paper, we will call a ring $R$ such that
$\Spec R$ has property FC an {\it FC-ring} or simply say {\it $R$ is
FC}. Such a condition on $R$ is useful when investigating questions
concerning heights of primes ideals. For example, using prime
avoidance one can show that if $P$ is a prime ideal of an FC-ring of
height at least $h$ then $P$ contains an ideal generated by $h$
elements which has height at least $h$.  And if $R$ is a quasi-local
FC-ring of dimension $d$ then every radical ideal is the radical of
a $d$-generated ideal.

An obvious question is whether the FC property passes to finitely
generated algebras.  Heinzer \cite{He} showed that if $R$ is an
FC-ring and $S$ is a finite $R$-algebra then $S$ is FC, while Ohm
and Pendleton prove that if $R$ is an FC-ring and $\dim R$ is finite
then $R[x]$ is FC. However, if $R$ is an infinite-dimensional
FC-ring $R[x]$ need not be FC.  This is illustrated by the following
example:

\begin{ex} \label{OPex} {\rm (\cite[Example 2.9]{OP})
Let $V$ be a valuation domain of countably infinite rank and let
$\Spec R=\{P_i\mid i\in \mathbb N_0\}$ where $P_{i}\subset
P_{i+1}$ for all $i$. Clearly, V is an FC-ring. For each $i\in
\mathbb N_0$ let $a_i\in P_{i+1}\setminus P_i$.  Let $x$ be an
indeterminate and set
$$f_i=a_i\prod_{j=0}^{i}(a_jx-1)\in V[x]$$
for each $i\in \mathbb N_0$.  Then the ideal $I=(\{f_i\mid i\in
\mathbb N_0\})V[x]$ has infinitely many minimal primes (namely,
$(P_i,a_ix-1)V[x]$ for $i\in \mathbb N_0$). Hence, $V[x]$ is not
an FC-ring.}
\end{ex}

In this paper we examine the weaker property  that every {\it
finitely generated} ideal of the ring has finitely many minimal
primes.  In many applications of FC (such as the ones mentioned in
the first paragraph), one only needs the FC property for finitely
generated ideals.   We show that this weaker property passes to
polynomial extensions (in any number of variables) as well as to
finitely presented algebras.

This paper arose in connection to the author's work on the
Cohen-Macaulay property for non-Noetherian rings \cite{HM}. It is an
open question whether the Cohen-Macaulay property (as defined in
\cite{HM}) passes to polynomial extensions.  While the answer in
general is probably `no', there is evidence to support a positive
answer if finitely generated ideals of the base ring have only
finitely many minimal primes.

\section{Main Results}

All rings in this paper are assumed to be commutative with identity.
For a ring $R$ and an $R$-module $M$ we let $\Supp M$ denote the
support of $M$ and $\Min M$ denote the minimal elements of $\Supp
M$.  We begin with the following definition:

\begin{df}{\rm Let $R$ be a ring. An $R$-module $M$ is said to be {\it FGFC} if for every
finitely generated $R$-submodule $N$ of $M$ it holds that $\Supp
M/N$ is the finite union of irreducible closed subsets of $\Spec R$;
equivalently, $\Min M/N$ is a finite set.  A ring $R$ is an {\it
FGFC-ring} if it is FGFC as an $R$-module.}
\end{df}

We make some elementary observations:

\begin{prop} Let $R$ be an FGFC ring.
\begin{enumerate}[(a)]
\item Every finitely presented $R$-module is FGFC. \item If $I$ is
a finitely generated ideal then $R/I$ is FGFC.  \item $R_S$ is an
FGFC ring for every multiplicatively closed set $S$ of $R$. \item
Every prime $p$ of finite height is the radical of a finitely
generated ideal. \item $R/p$ is an FGFC ring for every minimal
prime $p$ of $R$.
\end{enumerate}
\end{prop}

{\it Proof:}  Note that if $M$ is finitely presented then $\Min M=
\Min R/F_0(M)$, where $F_0(M)$ is the zeroth Fitting ideal of $M$
(\cite[Proposition 20.7]{E}). Thus, $\Min M$ is finite.  Since
every quotient of a finitely presented $R$-module by a finitely
generated submodule is also finitely presented, this proves (a).
Parts (b) and (c) are clear.  Part (d) is proved by induction and
prime avoidance. Part (e) follows from (b) and (d). \qed

\medskip

Clearly, every FC-ring is FGFC.  However, the converse is not
true.  The ring $V[x]$ of Example \ref{OPex} is not FC but is FGFC
by the following theorem:

\begin{thm} Let $R$ be a ring.  Then $R$ is an FGFC-ring if and only if
$R[x]$ is.
\end{thm}

{\it Proof:} Since $R\cong R[x]/(x)$ it is clear that $R$ is FGFC if
$R[x]$ is. Assume that $R$ is FGFC.  For a finitely generated ideal
$I$ of $R[x]$ let $d(I):=\min \{\sum_{i=1}^n\deg f_i\mid
I=(f_1,\dots,f_n)\}$.  If $d(I)=0$ then $I=JR[x]$ where $J$ is a
finitely generated ideal of $R$.   Hence $\operatorname{Min}_{R[x]}
R[x]/I=\{pR[x]\mid p\in \Min R/J\}$, which is a finite set as $R$ is
FGFC. Let $I$ be a finitely generated ideal of $R[x]$ with $d(I)>0$
and assume the theorem holds for all finitely generated ideals $J$
of $R[x]$ such that $d(J)<d(I)$. Let $f_1,\dots,f_n$ be a set of
generators for $I$ such that $\sum_{i=1}^n\deg f_i=d(I)$. Without
loss of generality we may assume $\deg f_n=\min_{1\le i\le n}\{\deg
f_i\mid \deg f_i>0\}$. Let $c$ be the leading coefficient of $f_n$
and $Q$ a minimal prime containing $I$. If $c\in Q$ then $Q$ is
minimal over $I'=(f_1,\dots,f_{n-1}, f_n-cx^{\deg f_n}, c)\supset
I$. As $d(I')<d(I)$ we have that $\operatorname{Min}_{R[x]} R[x]/I'$
is a finite set. Suppose that $c\not\in Q$.  Then $Q_c$ is minimal
over $I_c$.  We claim that $\operatorname{Min}_{R_c[x]}R_c[x]/I_c$
is finite.

\medskip

{\it Case 1:} $\deg f_i\ge \deg f_n$ for some $i<n$.

In this case, since $f_n$ is monic in $R_c[x]$ we can replace the
generator $f_i$ by one of smaller degree (using $f_n$).  Hence,
$d(I_c)<d(I)$ and $\operatorname{Min}_{R_c[x]}R_c[x]/I$ is a finite
set.

\medskip

{\it Case 2:} $\deg f_i<\deg f_n$ for all $i<n$.

By assumption on $\deg f_n$ we have $f_1,\cdots, f_{n-1}$ have
degree zero. By replacing $R$ with $R_c/(f_1,\dots,f_{n-1})R_c$
(which is still FGFC), we can reduce to the case $I=(f(x))$ where
$f(x)$ is a monic polynomial of positive degree.  As $S=R[x]/(f(x))$
is a free $R$-module, the going-down theorem holds between $S$ and
$R$. Therefore, every minimal prime of $S$ contracts to a minimal
prime of $R$. Further, since $S$ is finite as an $R$-module, there
are only finitely many primes of $S$ contracting to a given prime of
$R$. Hence, $\operatorname{Min}_{R[x]} R[x]/(f(x))$ is finite. \qed

\medskip

As a consequence, we get the following:

\begin{cor} Let $R$ be an FGFC-ring and $X$
a (possibly infinite) set of indeterminates over $R$. Then $R[X]$
is an FGFC-ring.
\end{cor}

{\it Proof:} Let $I=(f_1,\dots,f_n)$ be a finitely generated ideal
of $R[X]$.  Then there there exists $x_1,\dots,x_m\in X$ such that
$f_i\in R[x_1,\dots,x_m]$ for all $i$.  Let $S=R[x_1,\dots,x_m]$,
$J=(f_1,\dots,f_n)S$ and $X'=X\setminus \{x_1,\dots,x_m\}$.  By
the theorem (and induction) we have $\operatorname{Min}_SS/J$ is
finite. Furthermore, since $I=JS[X']$, every prime minimal over
$I$ is of the form $pS[X']$ for some $p\in
\operatorname{Min_S}S/J$. \qed

\medskip

It should be clear that, unlike FC, the FGFC property does not
pass to arbitrary finite ring extensions.  For example, if $V[x]$
and $I$ are as in Example \ref{OPex}, then $V[x]$ is FGFC but
$V[x]/I$ is not.  However, FGFC does pass to finitely presented
algebras:

\begin{cor} \label{fpalgebra} Let $R$ be an FGFC-ring.  Then any $R$-algebra which is
finitely presented as an $R$-module is an FGFC-ring.
\end{cor}

{\it Proof:} Let $S$ be an $R$-algebra which is finitely presented
over $R$.  Then certainly $S\cong R[x_1,\dots,x_n]/I$ for some
ideal $I$ of $R[x_1,\dots,x_n]$. By the theorem it is enough to
show that $I$ is finitely generated. Since $S$ is a finite
$R$-module, for each $i=1,\dots,n$ there exists a monic polynomial
$f_i(t)\in R[t]$ such that $f_i(x_i)\in I$. Let
$J=(f_1(x_1),\dots,f_n(x_n))R[x_1,\dots,x_n]$, $L=I/J$, and
$T=R[x_1,\dots,x_n]/J$.  It suffices to show that $L$ is a
finitely generated ideal of $T$.  In fact, since $T$ is a free
$R$-module of finite rank and $T/L$ is finitely presented over
$R$, it follows by the snake lemma that $L$ is finitely generated
as an $R$-module. \qed

\medskip

Given an integral ring extension $R\subset S$, it is clear that if
$S$ is FGFC then so is $R$.  For, any prime minimal over a finitely
generated ideal $I$ of $R$ is contracted from a prime of $S$ minimal
over $IS$.  However, FGFC does not in general ascend from $R$ to $S$
even if $R$ is coherent and $S$ is finite over $R$.  To see this,
one can again modify the example of Ohm and Pendleton:

\begin{ex} {\rm Let $R=V[x]$ and $I$ be as in Example \ref{OPex}.  Let
$S=R[y]/(yI,y^2-y)$.  Then $R$ is a coherent FGFC-domain (cf.
\cite{GV}), $R\subset S$ and $S$ is a finite $R$-module.  The set of
minimal primes of $S$ is $\{(P,y-1)S\mid P\in \Min R/I\}$, which is
infinite.  Hence, $S$ is not FGFC.}
\end{ex}

As a final result, we prove that FGFC does ascend from coherent
domains to finite torsion-free algebras:

\begin{prop} Suppose $R$ is
a coherent domain and $S$ a finite $R$-algebra which is
torsion-free as an $R$-module. If $R$ is FGFC then so is $S$.
\end{prop}

{\it Proof:} The hypotheses imply that $S$ is isomorphic as an
$R$-module to a finitely generated submodule of $R^n$ for some
$n$.  As $R$ is coherent, $S$ is finitely presented as an
$R$-module.  By Corollary \ref{fpalgebra}, $S$ is an FGFC ring.
\qed

\end{document}